\documentclass{amsart}
\usepackage{graphicx,latexsym,amssymb}

\newtheorem{thm}{Theorem}[section]
\newtheorem{lmm}{Lemma}[section]
\newtheorem{cor}{Corollary}[section]
\newtheorem{prop}[thm]{Proposition}
\newtheorem{defn}[thm]{Definition}
\newtheorem{rem}[thm]{Remark}
\newtheorem{example}[thm]{Example}

\newcommand{\C}{\mathbb{C}}
\newcommand{\PP}{\mathbb{P}}

\newcommand{\ben}{\begin{enumerate}}
\newcommand{\een}{\end{enumerate}}
\newcommand{\ble}{\begin{lem}}
\newcommand{\ele}{\end{lem}}
\newcommand{\bth}{\begin{thm}}
\renewcommand{\eth}{\end{thm}}
\newcommand{\bpr}{\begin{prop}}
\newcommand{\epr}{\end{prop}}
\newcommand{\bco}{\begin{cor}}
\newcommand{\eco}{\end{cor}}
\newcommand{\bde}{\begin{defn}}
\newcommand{\ede}{\end{defn}}
\newcommand{\brem}{\begin{rem}}
\newcommand{\erem}{\end{rem}}
\newcommand{\bexm}{\begin{example}}
\newcommand{\eexm}{\end{example}}

\begin{document}

\title{Rational transformations of systems of commuting nonselfadjoint operators}

\author{Alexander Shapiro and Victor Vinnikov}

\address{Alexander Shapiro,
Department of Mathematics and Statistics, Bar-Ilan University, Ramat-Gan 52900, Israel}
\email{sapial@math.biu.ac.il}

\address{Victor Vinnikov,
Department of Mathematics, Ben-Gurion University of the Negev, Beer-Sheva 84105, Israel}
\email{vinnikov@cs.bgu.ac.il}

\thanks{Partially supported by  EU-network HPRN-CT-2009-00099(EAGER), (The Emmy Noether Research Institute for Mathematics and the Minerva Foundation of Germany), the Israel Science Foundation grant \# 8008/02-3 (Excellency Center "Group  Theoretic Methods in the Study of Algebraic Varieties").}

\maketitle

\begin{abstract}
The work of M. S. Liv\v{s}ic and his collaborators in operator theory associates to a system of commuting
nonselfadjoint operators an algebraic curve. Guided by the notion of rational transformation of algebraic
curves, we define the notion of a rational transformation of a system of commuting nonselfadjoint
operators.
\end{abstract}

\section*{Introduction}

The work of M. S. Liv\v{s}ic and his collaborators in operator theory associates to a system of commuting
nonselfadjoint operators an algebraic curve (called the discriminant curve) given by a determinantal
representation, see \cite{LMKV}. This discovery leads to a very fruitful interplay between operator theory
and algebraic geometry: problems of operator theory lead to problems of algebraic geometry and vice versa.

A natural problem in operator theory is to define properly the notion of a rational transformation of a
system of commuting nonselfadjoint operators. This arises whenever one wants to study the algebra generated
by a given system of commuting nonselfadjoint operators. It may also allow representing the given system of
commuting nonselfadjoint operators in terms of another system which is simpler in some sense (e.g., it
contains fewer operators, or the operators have a smaller nonhermitian rank). A related problem in
algebraic geometry is to find an image of an algebraic curve given by a determinantal representation under
a rational transformation. To solve this problem we constructed elimination theory for pairs of polynomials
along an algebraic curve given by a determinantal representation. While our results are for polynomials in
two variables, plane algebraic curves, and pairs of operators, the generalization to polynomials in $d$
variables, algebraic curves in the $d$-dimensional space, and $d$-tuples of operators should be, for the
most part, relatively straightforward.

For the simplest case, when a system of operators consists of a single operator and the discriminant curve
is a line, these problems were solved by N.~Kravitsky using the classical elimination theory, see \cite{K}.
Our original objective was to find an analogue of the constructions of \cite{K} in the general case.

In the first section we formulate a problem of rational transformations of commuting nonselfadjoint
operators in a framework of commutative vessels, see \cite{LMKV}.

In the second section we discuss classical elimination theory and elimination theory along an algebraic
curve.

In the first part of Section 2 we remind classical notions and results of the elimination theory and
describe the image of a line under a rational transformation via a notion of Bezout matrix.

In the second part of Section 2 we introduce the notion of Vandermonde vector on a plane algebraic curve
and generalize the notion of Bezout matrix. We formulate analogues of the results of the classical
elimination theory for the elimination theory along a plane algebraic curve. We describe an image of a
plane algebraic curve given by a determinantal representation under a rational transformation via the
notion of generalized Bezout matrix.

In Section 3 we solve a problem formulated in Section 1 using elimination theory along a plane algebraic
curve.
\section{Statement of the problems}

In the recent years a significant progress has been made in the study of commuting non-selfadjoint
operators, see \cite{LMKV}. It turned out that associated to a system of commuting non-selfadjoint
operators with finite--dimensional imaginary part there is a certain algebraic variety called the
discriminant variety. Guided by the standard notion of rational transformation of algebraic varieties it
thus becomes a natural problem to study rational transformations of systems of commuting non-selfadjoint
operators, which is the goal of this section.

To formulate the main problem more precisely, we have to introduce some notation. We shall use the
framework of commutative vessels \cite{LMKV} which turns out to be very convenient in the study of
commuting non-selfadjoint operators; it generalizes to the multi--operator case the framework of
colligations (nodes) which has been extensively used in the study of a single non-selfadjoint (or
non-unitary) operator, see, e.g., \cite{B}.

Let $H$ be a Hilbert space (finite-- or infinite--dimensional) and let $E$ be a finite--dimensional Hilbert
space.

\begin{defn}
A commutative vessel is a collection
\begin{eqnarray}
V=(A_{i}, H, {\Phi}, E, {\sigma}_{i}, {\gamma}^{in}_{ij}, {\gamma}^{out}_{ij}) \nonumber
\end{eqnarray}
where $i{\,}, j=1{\,}{\dots}{\,}k{\,},{\,} A_{i}:H{\,}{\to}{\,}H$ are {\,}bounded {\,}linear {\,}
commuting {\,}operators;\\
${\Phi}${\,} is{\,} a{\,} bounded{\,} linear{\,} mapping {\,}from {\,}$H${\,} to{\,} $E$
{\,}with{\,} the{\,} adjoint{\,} mapping\\
${\Phi}^{*}:E{\to}H; {\,}{\sigma}_{i}, {\gamma}^{in}_{ij}, {\gamma}^{out}_{ij}$ are bounded selfadjoint
operators in $E$, such that\\ ${\gamma}^{in}_{ij}=-{\gamma}^{in}_{ji},
{\gamma}^{out}_{ij}=-{\gamma}^{out}_{ji}$, and
\begin{eqnarray}
{\Phi}^*{\sigma}_{i}{\Phi} & = & \frac{ 1 }{ i }(A_{i}-A^*_{i}) \nonumber \\
{\gamma}^{in}_{ij}{\Phi} & = & {\sigma}_i{\Phi}A^*_{j}-{\sigma}_{j}{\Phi}A^*_{i} \nonumber \\
{\gamma}^{out}_{ij}{\Phi} & = & {\sigma}_i{\Phi}A_{j}-{\sigma}_{j}{\Phi}A_{i}\nonumber \\
{\gamma}^{out}_{ij} & = & {\gamma}^{in}_{ij}+i({\sigma}_{i}{\Phi}{\Phi}^*{\sigma}_{j}-
{\sigma}_{j}{\Phi}{\Phi}^*{\sigma}_{i}) \nonumber
\end{eqnarray}
\end{defn}

The first problem is to define properly rational transformations of the vessel. Let $r_{1}
(y_{1},{\dots},y_{k}),{\dots},r_{l}(y_{1},{\dots},y_{k}){\,}$ be $l$ real rational functions of $k$
variables such that $A'_{1}=r_{1}(A_{1},{\dots},A_{k}),{\dots}, A'_{l}=r_{l}(A_{1},{\dots},A_{k})$ are
defined.

{\bf Problem 1} Given $V$ as above, find ${\Phi}', E' ,{\sigma}'_{i}, {\gamma}'^{in}_{ij},
 {\gamma}'^{out}_{ij}$ so that $V'=(A'_{i},H,{\Phi}',\\ E' ,{\sigma}'_{i}, {\gamma}'^{in}_{ij},
 {\gamma}'^{out}_{ij})$ is a vessel.

For the case $k=1$ this problem has been solved by Naftali Kravitsky (see \cite{LMKV}, \cite{K}) using the
notion of the Bezout matrix of two polynomials.

As we have mentioned in the beginning rational transformations of vessels should be connected with rational
transformations of the associated discriminant varieties. We proceed to state the corresponding problems
restricting ourselves for simplicity to the case of two commuting non-selfadjoint operators ($k=2$). We
will denote ${\gamma}^{in}_{12}={\gamma}^{in}, {\gamma}^{out}_{12}={\gamma}^{out}$. We thus consider a two
operator vessel
\begin{displaymath}
V=(A_{1},A_{2}, H, {\Phi}, E, {\sigma}_{1}, {\sigma}_{2}, {\gamma}^{in}, {\gamma}^{out}).
\end{displaymath}
Next theorem was proved in \cite{LMKV}.
\begin{thm}
\begin{displaymath}
\det(y_{1}{\sigma}_{2}-y_{2}{\sigma}_{1}+{\gamma}^{in})
=\det(y_{1}{\sigma}_{2}-y_{2}{\sigma}_{1}+{\gamma}^{out})
\end{displaymath}
\end{thm}
\begin{defn} The polynomial ${\Delta}(y_{1},y_{2})$ defined by
\begin{displaymath}
{\Delta}(y_{1},y_{2})=\det(y_{1}{\sigma}_{2}-y_{2}{\sigma}_{1}+{\gamma}^{in})
=\det(y_{1}{\sigma}_{2}-y_{2}{\sigma}_{1}+{\gamma}^{out})
\end{displaymath}
is called the discriminant polynomial of the vessel.
\end{defn}

The discriminant polynomial is an important invariant of the corresponding vessel. In particular we have
the generalized Cayley-Hamilton Theorem: $${\Delta}(A_{1},A_{2})=0$$ (on the so-called principal subspace
of the vessel).

We shall assume for simplicity that ${\Delta}$ is an irreducible polynomial and denote by $C$ the
irreducible plane algebraic curve with the equation ${\Delta}=0$. $C$ is called the discriminant curve. It
follows from the irreducibility of ${\Delta}$ that if we shall define for every $(y_{1},y_{2})$ on $C$ the
following two subspaces of $E$
\begin{eqnarray}
E^{in}(y_{1},y_{2}) & = & ker(y_{1}{\sigma}_{2}-y_{2}{\sigma}_{1}+
{\gamma}^{in}) \nonumber \\
E^{out}(y_{1},y_{2}) & = & ker(y_{1}{\sigma}_{2}-y_{2}{\sigma}_{1}+
{\gamma}^{out}) \nonumber
\end{eqnarray}
then $E^{in}$ and $E^{out}$ are line bundles on the Riemann surface $X$ of the algebraic curve $C$ .

\begin{defn} Let $C$ and $C'$ be irreducible plane algebraic curves. $C$ and $C'$ are called
birationally isomorphic if there exist rational mappings of $C$ into $C'$ and of $C'$ into $C$ that are
inverse to one another. Equivalently there exists a rational mapping from $C$ to $C'$ which is one-to-one
almost everywhere (i.e., except at a finite number of points).
\end{defn}

Let $r_{1},r_{2}$ be real rational functions, such that $r=(r_{1},r_{2})$ is a birational isomorphism of an
irreducible plane algebraic curve $C$ to an irreducible plane algebraic curve $C'$. Let $V$ be a vessel
with discriminant curve $C$ and let $V'$ be the vessel constructed from it in Problem 1.

{\bf Problem 2}

{\bf 1.} Prove that the discriminant polynomial of $V'$ is irreducible and that the discriminant curve of
$V'$ is $C'$.

{\bf 2.} Prove that $E'^{in}$ is isomorphic to $E^{in}$ and $E'^{out}$ is isomorphic to $E^{out}$ (more
precisely isomorphic up to a certain twist; note that since $C$ and $C'$ are birationally isomorphic they
have the same associated Riemann surface).

Technical remark: If $C$ or $C'$ have singularities than in Problem 2 we ought to consider only  vessels
whose corresponding determinant representations are maximal, that is the dimension of $E^{in}$ and
$E^{out}$ is maximal for all singular points of the curve, see \cite{BV}.

One motivation for Problem 2 comes from the case $C={\PP}^{1}$, that is $C$ is a straight line. $C'$ is
then a rational plane curve and $(r_{1},r_{2})$ is its Luroth (i. e. almost everywhere one-to-one) rational
parametrization. In this case all the statements of Problem 2 were actually proved by Naftali Kravitsky,
see \cite{K}. An additional evidence is provided by the theory of functional models where one constructs at
once the model for any operator in the ring generated by $A_{1}$ and $A_{2}$ , see \cite{AV}.
\section{Classical elimination theory and elimination theory along a curve}
\subsection{Classical elimination theory}
\subsubsection{Vandermonde vectors}
Let us recall the main goal of elimination theory. Given $n+1$ (nonhomogeneous) polynomials in $n$
variables we want to find necessary and sufficient conditions (in terms of the coefficients) for these
polynomials to have a common zero (and furthermore to determine the number of common zeroes, counting
multiplicities, if they exist), see \cite{M}.

In the classical case we consider two (nonhomogeneous) polynomials in one variable. For such a pair of
polynomials one may construct so called Bezout matrix. Entries of this matrix depend only on coefficients
of the polynomials. Determinants of this matrices equal to zero if and only if the polynomials have a
common zero.

We will consider classical elimination theory in a framework of Vandermonde vectors. This approach allowed
us to give new proofs of classical results of elimination theory for two polynomials in one variable. The
advantage of this approach is that these proofs can be generalized for elimination theory for two
polynomials in two variables along a plane algebraic curve, see \cite{SV}.

Let $V_n(x)$ be a Vandermonde vector of the length $n$:
$$V_n(x) =
\left(
\begin{array}{c}
1 \\
x \\
\dots \\
x^{n-1}
\end{array}
\right) =
\left( x^i \right)_{i=0}^{n-1}$$

\begin{thm}
If $x_1, x_2, \dots , x_n$ are pairwise distinct then vectors\\ $V_n(x_1), V_n(x_2), \dots , V_n(x_n)$ are
linearly independent.
\end{thm}
\subsubsection{Bezout matrices}
\begin{lmm}
For every two polynomials in one variable $p(x)$ and $q(x)$ of degree $n$ there exists uniquely determined
$n \times n$ symmetric matrix $B(p,q) = \left( b_{ij} \right)_{i,j=0}^n$ such that $p(x)q(y) - q(x)p(y) =
\sum b_{ij}x^i(x-y)y^j$.
\end{lmm}
\begin{cor}
$p(x)q(y) - q(x)p(y) = V^T_n(x)(x-y)B(p,q)V_n(y)$.
\end{cor}
\begin{thm}
The dimension of the kernel of Bezout matrix of two polynomials $p$ and $q$ of degree $n$ is equal to the
number of common zeroes of these polynomials and $|\det S(p,q)|=|\det B(p,q)|$.
\end{thm}
The determinant of the Bezout matrix is called a bezoutian.
\subsubsection{Rational transformations}
We will use the notion of Bezout matrix to describe an image of ${\C}$ under a rational transformation. Let
us consider three polynomial in one variable $p_0(x)$, $p_1(x)$ and $p_2(x)$. These polynomials define a
rational transformation $r: {\C} \to {\C}^2$: by $r(x)=(\frac{p_1(x)}{p_0(x)}, \frac{p_2(x)}{p_0(x)})$.

Next theorem is due to Kravitsky, \cite{K}.
\begin{thm}
The image of ${\C}$ under the rational transformation $r$ is a curve given by the next determinantal
representation:
$$\det(B(p_0,p_1)x + B(p_0,p_2)y + B(p_1,p_2)) = 0$$
\end{thm}
\subsection{Elimination theory along a curve}
\subsubsection{Vandermonde vectors on a curve}
It is an unfortunate fact that the easy, useful and beautiful constructions of the single variable case do
not generalize in any straightforward fashion to the case of two or more variables. Some generalizations of
the classical results using much more sophisticated algebraic techniques have been obtained in recent years
by Gelfand, Kapranov, Zelevinsky, see \cite{GKZ}, and Jouanolou, see \cite{J}. In \cite{SV} we followed a
different, somewhat ``asymmetrical'' approach. Namely, we chosed one of our three polynomials and viewed it
as defining a plane algebraic curve. Our goal became to generalized the notion of Bezout matrix for three
(homogeneous) polynomials in two variables on the curve defined by the third polynomial.

To solve this problem it turns out to be essential to consider a plane algebraic curve given by a
determinantal representation
$${\Delta}(x_0,x_1,x_2)=\det(x_0D_0+x_1D_1-x_2D_2)$$
where $D_0$, $D_1$ and $D_2$ are $m \times m$ complex hermitian matrices.

The determinantal representation of a curve provides us with an additional structure. In the case of a
plane irreducible curve for every point of the curve $x=(x_0,x_1,x_2)$ there is a nontrivial subspace
$\ker(x_0D_0 + x_1D_1 - x_2D_2)$ (one--dimensional, because of the irreducibility of
${\Delta}(x_0,x_1,x_2)$, except at the singular points). We choose a vector $e$ in this subspace. This
additional structure allows us to define an appropriate analogue of Vandermonde vectors. Vandermonde
vectors on a curve (equipped with a determinantal representation) is given by
$$V_n(x,e) =
\left(
\begin{array}{c}
x_0^{n-1}e \\
x_0^{n-2}x_1e {\quad} x_0^{n-2}x_2e \\
\dots {\quad} \dots {\quad} \dots {\quad} \dots {\quad} \dots \\
x_0x_1^{n-2}e {\quad} \dots {\quad} \dots {\quad} \dots {\quad} \dots {\quad} x_0x_2^{n-2}e \\
x_1^{n-1}e {\quad} x_1^{n-2}x_2e {\quad} \dots {\quad} \dots {\quad} \dots {\quad} x_1x_2^{n-2}e
{\quad} x_2^{n-1}e
\end{array}
\right) =$$
$$=\left( x_0^{i_0}x_1^{i_1} x_2^{i_2} e \right)_{i_0+i_1+i_2=n-1} = \left( x^i e
\right)_{|i|=n-1}$$

It is clear that a Vandermonde vector $V_n(x,e)$ on a curve belongs to the space
${\C}^{m\frac{n(n+1)}{2}}$, where $m$ is the degree of ${\Delta}(x_0,x_1,x_2)$. We will denote this space
by $W_n$ and will call it the blown up space. Unlike the classical case all the Vandermonde vectors on a
curve generate not the whole blown space $W_n$ but only a certain subspace of $W_n$. Let us consider a
subspace
$$V_n = \{ (w_{j_0,j_1,j_2}) \in W_n :
D_0w_{j_0+1,j_1,j_2} + D_1w_{j_0,j_1+1,j_2} - D_2w_{j_0,j_1,j_2+1}=0 \}$$
where $j_0+j_1+j_2=n-1$. We will call $V_n$ the principal subspace.
It is clear that all Vandermonde vectors belong to $V_n$, that is
$V_n(x) \in V_n$ for every point $x=(x_0,x_1,x_2)$.
\begin{thm} The dimension of the principal subspace $V_n$ equals to $nm$.
\end{thm}
\begin{thm}
Vandermonde vectors in zeroes of any polynomial on a curve are linearly independent and generate the
principal subspace $V_n$.
\end{thm}

The blown space and the principal subspace play a key role in the elimination theory on a plane algebraic
curve. Namely, all the constructions of the classical elimination theory can be formally generalized in
terms of the blown spaces. After this formal generalization we restrict the new ``blown up'' constructions
to the corresponding principal subspaces and obtain complete analogues of the classical results.
\subsubsection{Generalized Bezout matrices}
Our goal in this paragraph is to define properly the generalized Bezout matrix for a pair of homogeneous
polynomials in three variables on a curve given by a determinantal representation.

As in the previous paragraph we consider the blown space $W_n$ and the principal subspace $V_n$. For a pair
of polynomials we will determine formal analogues of classical Bezout matrix in $W_n$ and then will
restrict this analogue on $V_n$.
\begin{lmm}
For every two homogeneous polynomials in three variables \\
$p(x_0,x_1,x_2)$ and $q(x_0,x_1,x_2)$ of degree $n$ there exist three
$\frac{n(n+1)}{2} \times \frac{n(n+1)}{2}$ symmetric matrices
${\beta}^{10} = (b^{10}_{ij})$, ${\beta}^{20} = (b^{20}_{ij})$
and ${\beta}^{12} = (b^{12}_{ij})$ such that
$p(x_0,x_1,x_2)q(y_0,y_1,y_2) - q(x_0,x_1,x_2)p(y_0,y_1,y_2) =$
$$\sum_{|i|,|j|=n}
b^{10}_{ij} x^i (x_1y_0-x_0y_1) y^j +
b^{20}_{ij} x^i (x_2y_0-x_0y_2) y^j +
b^{12}_{ij} x^i (x_1y_2-x_2y_1) y^j $$
\end{lmm}

Let us define a $\frac{n(n+1)}{2}m \times \frac{n(n+1)}{2}m$
matrix $B(p,q)$ on the space $W_n$: \\
$B(p,q) = {\beta}^{12} \otimes D_0 + {\beta}^{10} \otimes D_1 + {\beta}^{20} \otimes D_2$.

Let us consider the restriction of $B(p,q)$ on the subspace generated by Vandermonde vectors on the curve:
$$B'(p,q) = \mathcal{P}_{V_n} B(p,q) \mathcal{P}_{V_n}$$

We will call $B'(p,q)$ Bezout matrix of polynomials $p$ and $q$ along the algebraic curve given by the
determinantal representation ${\Delta}$.

Next useful lemma follows from the definition of the generalized Bezout matrix.
\begin{lmm}If
$V(x)=(x^i e)=(x_0^{i_0}x_1^{i_1}x_2^{i_2} e)$ and $V(y)=(y^i h)=(y_0^{i_0}y_1^{i_1}y_2^{i_2} h)$
are two Vandermonde vectors on a curve then\\
$V^T(y)B'(p,q)V(x)=(p(x)q(y)-q(x)p(y)) h^T e$
\end{lmm}
\begin{thm}
The dimension of the kernel of Bezout matrix of two polynomials along the curve is equal to the number of
common zeroes of these polynomials along the curve.
\end{thm}
\begin{cor}
Two polynomials in two variables have common zero on a curve if and only if the determinant of their Bezout
matrix along this curve is equal to zero.
\end{cor}
\subsubsection{Rational transformation of a curve}
As was mentioned in the paragraph 2.1.3 the image of the strait line under a rational transformation can be
described via the notion of the Bezout matrix. In this paragraph we consider a rational transformation of a
curve given by a determinantal representation and describe the image of a curve via the notion of the
generalized Bezout matrix.

As before, we start from a real plane algebraic curve $C$ defined by a homogeneous polynomial in three
variables ${\Delta}(x_0,x_1,x_2)$:
$$C = \{ (x_0,x_1,x_2) \in {\PP}^2: {\Delta}(x_0,x_1,x_2)=0 \}$$
and a determinantal representation of this curve:
$${\Delta}(x_0,x_1,x_2)= \det (x_0D_0 + x_1D_1 - x_2D_2)$$

Let us consider three homogeneous polynomials in three variables $p_0$, $p_1$, $p_2$
of degree $n$ and an image $r(C)$ of the curve $C$ under the rational transformation
$r=(p_0,p_1,p_2)$:
$$r(C) = \{ (p_0(x),p_1(x),p_2(x)) \in {\PP}^2: x \in C\}$$

To find a determinantal representation of the curve $r(C)$ we consider three
generalized Bezout matrices: $B'(p_0,p_1)$, $B'(p_0,p_2)$ and $B'(p_1,p_2)$.
\begin{lmm}
The curve given by determinantal representation\\
$\det (x_0B'(p_1,p_2) + x_1B'(p_0,p_2) - x_2B'(p_0,p_1))$ contains $r(C)$.
\end{lmm}
But if a basepoint $y=(y_0,y_1,y_2)$ of the rational transformation $r=(p_0,p_1,p_2)$
belongs to the curve $C$, that is if there exists a common zero of polynomials $p_0$, $p_1$ and $
p_2$
on the curve, then $\det (x_0B'(p_1,p_2) + x_1B'(p_0,p_2) - x_2B'(p_0,p_1))=0$
because $B'(p_i,p_j)V(y)=0$ for every $i$ and $j$.

Therefore we consider a subspace $\bar V_n$ of the principal subspace $V_n$ which is perpendicular to the
Vandermonde vectors in basepoints and restrict generalized Bezout matrices on this subspace:\\
$\bar V_n = <V_n(y): p_0(y)=p_1(y)=p_2(y)=0>^{\perp}$, $\bar B(p_i,p_j) = \mathcal{P}_{\bar V_n}
B'(p_i,p_j) \mathcal{P}_{\bar V_n}$.
\begin{thm}
The curve $r(C)$ has the determinantal representation\\
$\det (x_0 \bar B(p_1,p_2) + x_1 \bar B(p_0,p_2) - x_2 \bar B(p_0,p_1))$.
\end{thm}
\section{Solution of the problems}

As before, we consider a two operator vessel
$$V=(A_1,A_2, H, {\Phi}, E, {\sigma}_1, {\sigma}_2, {\gamma}^{in},
{\gamma}^{out})$$
with the discriminant curve $C$ defined by the determinantal representation
$${\Delta}(y_1,y_2)=\det(y_1{\sigma}_2-y_2{\sigma}_1+{\gamma}^{in})$$
and a rational transformation
$$r_1(y_1,y_2) =\frac{p_1(y_1,y_2)}{p_0(y_1,y_2)}{\quad} ,{\quad}
r_2(y_1,y_2) =\frac{p_2(y_1,y_2)}{p_0(y_1,y_2)}$$

Let us define by $n$ the maximal degree of the polynomials
$p_0(y_1,y_2)$, $p_1(y_1,y_2)$ and $p_2(y_1,y_2)$.
Let us homogenize all three polynomials up to the degree $n$. Thus we consider
three homogeneous polynomials in three variables of degree $n$:
$P_0(x_0,x_1,x_2)$, $P_1(x_0,x_1,x_2)$ and $P_2(x_0,x_1,x_2)$.

As in the second part of Section 2 we consider the blown space
$W_n={\bf \rm C}^{m\frac{n(n+1)}{2}}$,
the principal subspace $V_n \subset W_n$ defined by formula
$$V_n=\{ (w_{i_1,i_2}) \in W_n:
{\sigma}_1w_{i_1+1,i_2}-{\sigma}_2w_{i_1,i_2+1}+{\gamma}^{in}w_{i_1,i_2}=0 \}$$
and Bezout matrices of the polynomials $P_j$ and $P_k$ along the curve $C$ denoted
by $B'(P_j,P_k)$.

Let us define
\begin{eqnarray}
E' & = & V_n \nonumber \\
{\Phi}' & = & P_{E'}({\Phi}A^{*i_1}_1A^{*i_2}_2
p^{-1}_0(A^*_1A^*_2))_{i_1i_2} \nonumber \\
{\sigma}'_{1} & = & B'(P_0,P_1) \nonumber \\
{\sigma}'_{2} & = & B'(P_0,P_2) \nonumber \\
{\gamma}'^{in} & = &  B'(P_1,P_2) \nonumber \\
{\gamma}'^{out} & = & {\gamma}'^{in}+i({\sigma}'_1{\Phi}'{\Phi}'^*{\sigma}'_2-
{\sigma}'_2{\Phi}'{\Phi}'^*{\sigma}'_1) \nonumber \\
V' & = & (A'_1, A'_2, H, {\Phi}', E', {\sigma}'_1, {\sigma}'_2,
{\gamma}'^{in}, {\gamma}'^{out}) \nonumber
\end{eqnarray}
\begin{thm}
$V'$ is a vessel.
\end{thm}
{\bf Proof } If the statement of the theorem is true for the transformation defined by polynomials $p_0,
p_1, p_2$ and for the transformation defined by polynomials $p_0, q_1, q_2$ then the statement is also true
for the transformation defined by polynomials $p_0, p_1+q_1, p_2+q_2$. Hence, it is sufficient to prove the
theorem for the transformation defined by pairs of monomials $p_0, x_1^{i_1}x_2^{i_2}, x_1^{j_1}x_2^{j_2}$.
For monomials this theorem follows recursively from the identity\\
$x_1^ny_2^m - x_2^my_1^n =
x_1(x_1^{n-1}y_2^m - x_2^my_1^{n-1}) +
(x_1^{n-1}y_2^m - x_2^my_1^{n-1})y_1 + \\
x_1(x_1^{n-2}y_2^m - x_2^my_1^{n-2})y_1$
and the fact that operators\\
$T_1:E'_{n-1} \to E'_n$, $e_{i_1,i_2} \mapsto e_{i_1+1,i_2}$ and $T_2:E'_{n-1} \to E'_n$, $e_{i_1,i_2}
\mapsto e_{i_1,i_2+1}$ preserve the conditions of the vessel. Theorem is proved.
\begin{thm}
The discriminant variety of $V'$ contains $C'$.
\end{thm}
{\bf Proof } Let us consider the point ${\lambda} = ({\lambda}_1,{\lambda}_2)$
of the discriminant curve $C$.\\
By definition of the discriminant polynomial
$\det({\lambda}_1{\sigma}_2-{\lambda}_2{\sigma}_1+{\gamma}^{in}) = 0$.
Hence, there exist a non-zero vector ${\phi}$ such that
$({\lambda}_1{\sigma}_2-{\lambda}_2{\sigma}_1+{\gamma}^{in}){\phi} = 0$.\\
Let us consider the birational transformation $r = (r_1, r_2)$ of a vessel\\
${\Lambda} = ({\lambda}_1,{\lambda}_2, H, {\phi}, {\C}^n, {\sigma}_1, {\sigma}_2, {\gamma}^{in},
{\gamma}^{out})$ as above. It follows from\\the Theorem 3.1 that ${\Lambda}' = ({\lambda}'_1, {\lambda}'_2,
H, {\phi}', {\C}^{n^2}, {\sigma}'_1,
{\sigma}'_2, {\gamma}'^{in}, {\gamma}'^{out})$ is a vessel.\\
By definition of a vessel
${\sigma}'_2{\phi}'{\lambda}'_1 - {\sigma}'_1{\phi}'{\lambda}'_2
+ {\gamma}^{in}{\phi}' = 0$. Therefore\\
$\det({\lambda}'_1{\sigma}'_2 - {\lambda}'_2{\sigma}'_1
+ {\gamma}^{in} = 0$ and the point ${\lambda}' = (r_1({\lambda}),r_2({\lambda}))$
belongs to the discriminant variety of $V'$. The theorem is proved.

We obtained a solution of problem 1. But if a basepoint of the rational transformation $r_1$, $r_2$ belongs
to the discriminant curve of the vessel $V$ then
the discriminant variety of the vessel $V'$ is the whole plane. That is if there \\
exist a point ${\lambda}=({\lambda}_1,{\lambda}_2)$
such that $\det({\lambda}_2{\sigma}_1 - {\lambda}_1{\sigma}_2 + {\gamma})=0$ and \\
$p_i({\lambda})=0 , i=0,1,2$, then for every $y = (y_1 , y_2)$
$\det(y_2{\sigma}'_1 - y_1{\sigma}'_2 + {\gamma}')=0$.

Therefore, as in Section 2, part 2, paragraph 3, we consider
\begin{eqnarray}
E'' & = & <({\lambda}_{1}^{i_{1}}{\lambda}_{2}^{i_{2}}e)_{i_{1}i_{2}}
:{p_{n}({\lambda})=0},
({\lambda}_{2}{\sigma}_{1} - {\lambda}_{1}{\sigma}_{2} + {\gamma})e=0>^{\perp} =
\bar V_n \nonumber \\
{\Phi}'' & = &  \mathcal{P}_{E''}{\Phi}' \nonumber \\
{\sigma}''_{1} & = & \mathcal{P}_{E''}{\sigma}'_{1} \mathcal{P}_{E''} = \bar B(P_0,P_1) \nonumber \\
{\sigma}''_{2} & = & \mathcal{P}_{E''}{\sigma}'_{2} \mathcal{P}_{E''} = \bar B(P_0,P_2) \nonumber \\
{\gamma}''^{in} & = & \mathcal{P}_{E''}{\gamma}'^{in} \mathcal{P}_{E''} = \bar B(P_1,P_2) \nonumber \\
{\gamma}''^{out} & = & \mathcal{P}_{E''}{\gamma}'^{out} \mathcal{P}_{E''} \nonumber \\
V'' & = & (A'_{1},A'_{2}, H, {\Phi}'', E'', {\sigma}''_{1}, {\sigma}''_{2},
{\gamma}''^{in}, {\gamma}''^{out}) \nonumber
\end{eqnarray}
Next theorem follows from lemma 2.4.
\begin{thm}
$V''$ is a vessel and the discriminant variety of $V''$
contains $C'$.
\end{thm}

It follows from corollary 2.2 that $\det{\sigma}'_i=0$ if and only if the polynomials $p_0$ and $p_i$ have
a common zero on the curve $C$ and that $\det{\gamma}'=0$ if and only if the polynomials $p_1$ and $p_2$
have a common zero on the curve $C$. This fact implies next two theorems.

\begin{thm}
Discriminant variety of $V''$ is $C'$.
\end{thm}
\begin{thm}
$E^{in}$ is isomorphic to $E''^{in}$ and $E^{out}$ is isomorphic to $E''^{out}$. The~isomorphism maps the
vectors $e \in E^{in}({\lambda}_1, {\lambda}_2)$ to the vectors $({\lambda}^{i_1}_1 {\lambda}^{i_2}_2 e)_{0
\leq i_1+i_2 \leq n} \in E''^{in}$ and the vectors $e \in E^{out}({\lambda}_1, {\lambda}_2)$ to the vectors
$({\lambda}^{i_1}_1 {\lambda}^{i_2}_2 e)_{0 \leq i_1+i_2 \leq n} \in E''^{out}$.
\end{thm}


\clearpage

\begin{thebibliography}{}

\bibitem{AV} D. Alpay, V. Vinnikov, {\it Analogues d'espaces de de Branges sur des surfaces de Riemann},
C. R. Acad. Sci. Paris, t. 318, Serie I, p. 1077-1082, 1994.

\bibitem{B} M. S. Brodskii, {\it Triangular and Jordan representations of linear operators},
Transl. Math. Monographs 32, Amer. Math. Soc., Providence, 1970.

\bibitem{BV} Joseph A. Ball, Victor Vinnikov, {\it Zero-pole interpolation for matrix
meromorphic functions on an algebraic curve and transfer functions of 2d systems}, Acta Applicandae
Mathematicae. To appear.

\bibitem{GKZ} I. M. Gelfand, M. M. Kapranov, A. V. Zelevnsky, {\it Discriminants, resultants
and multidimensional determinants}, Boston ,Birkhauser, 1994.

\bibitem{J} J.--P. Jouanolou {\it Le formalisme du r\'{e}sultant}, Adv. Math. 90 (1991), p. 117-263

\bibitem{K} N. Kravitsky {\it On the discriminant function of two commuting non-selfadjoint
operators}, Integral Equations Operator Theory 3/1, p. 97-124, 1980.

\bibitem{LMKV} M. S. Liv\v{s}ic, A. S. Markus, N. Kravitsky, V. Vinnikov, {\it Theory of Commuting
Nonselfadjoint Operators}, Kluwer, Dordrecht, 1995.

\bibitem{M} F. S. Macaulay {\it The algebraic theory of modular systems}, New York,
Stechert--Hafner Service Agency, 1964.

\bibitem{SV} Shapiro A., Vinnikov V., {\it Rational transformations of algebraic curves and elimination theory},
Linear Algebra and its Applications, preprint.

\bibitem{V} V. Vinnikov, {\it Self--adjoint determinantal representations of real plane curves},
Math. Ann., 296 (1993), p. 453-473.
\end{thebibliography}
\end{document}